\RequirePackage{fix-cm}

\documentclass[smallextended]{svjour3}       
\smartqed  
\usepackage{latexsym,amsmath,amssymb}
\usepackage[english]{babel}

\newcommand{\ds}{\displaystyle}
\newcommand{\llr}{\Longleftrightarrow}
\newcommand{\lr}{\Longrightarrow}
\newcommand{\al}{\alpha}
\newcommand{\la}{\lambda}
\newcommand{\si}{\sigma}
\journalname{Journal of Theoretical Probability}
\begin{document}

\title{An orthogonal-polynomial approach to first-hitting times of birth-death
processes}

\titlerunning{first-hitting times of birth-death processes}    

\author{Erik A. van Doorn}

\institute{E.A. van Doorn \at
	Department of Applied Mathematics, University of Twente\\
	P.O. Box 217, 7500 AE Enschede, The Netherlands\\             
	\email{e.a.vandoorn@utwente.nl}
}

\date{Received: date / Accepted: date}

\maketitle

\begin{abstract}
In a recent paper in this journal Gong, Mao and Zhang, using the theory of
Dirichlet forms, extended Karlin and McGregor's classical results on
first-hitting times of a birth-death process on the nonnegative integers by
establishing a representation for the Laplace transform $\mathbb{E}[e^{sT_{ij}}]$
of the first-hitting time $T_{ij}$ for {\em any} pair of states $i$ and $j$, as
well as asymptotics for $\mathbb{E}[e^{sT_{ij}}]$ when either $i$ or $j$ tends
to infinity. It will be shown here that these results may also be obtained by
employing tools from the orthogonal-polynomial toolbox used by Karlin and
McGregor, in particular {\em associated polynomials} and {\em Markov's Theorem}.
\end{abstract}

\section{Introduction}

A {\em birth-death process\/} is a continuous-time Markov chain $\mathcal{X}
:= \{X(t),~t \geq 0\}$ taking values in $S := \{0,1,2,\ldots\}$ with $q$-matrix
$Q := (q_{ij},~i,j \in S)$ given by
\[
\begin{array}{@{}l}
q_{i,i+1} = \la_i,~~q_{i+1,i} = \mu_{i+1},~~q_{ii} = -(\la_i + \mu_i), \\
q_{ij} = 0,\quad |i-j|>1,
\end{array}
\]
where $\la_i > 0$ for $i \geq 0,~ \mu_i > 0$ for $i \geq 1$ and $\mu_0 \geq 0$.
Positivity of $\mu_0$ entails that the process may evanesce by escaping from
$S$, via state 0, to an absorbing state $-1$. Throughout this paper we will
assume that the transition probabilities
\[
P_{ij}(t) := \mathbb{P}(X(t)=j\,|\,X(0)=i), \quad i,j \in S,
\]
satisfy both the backward and forward Kolmogorov equations, and mostly also 
that they are uniquely determined by the {\em birth rates\/} $\la_i$ and
{\em death rates\/} $\mu_i$. Karlin and McGregor \cite{K57a} have shown
that the latter is equivalent to assuming 
\begin{equation}
\label{unique}
\sum_{n=0}^\infty \left(\pi_n + \frac{1}{\la_n\pi_n}\right) = \infty,
\end{equation}
where the $\pi_n$ are constants given by
\[
\pi_0 := 1 \mbox{~~and~~} \pi_n := 
	\frac{\la_0\la_1\ldots\la_{n-1}}{\mu_1\mu_2\ldots\mu_n}, ~~n > 0.
\]
We note that condition \eqref{unique} does not exclude the possibility of
{\it explosion}, escape from $S$, via all states larger than the
initial state, to an absorbing state $\infty$.

We denote by $T_{ij}$ the (possibly defective) first hitting time of state $j$,
starting in state $i\neq j$. Then, writing
\[
\hat{P}_{ij}(s) := \int_0^\infty e^{st}P_{ij}(t)dt, \quad s<0, 
\]
and
\[
\hat{F}_{ij}(s) := \mathbb{E}[e^{sT_{ij}}] =
\int_0^\infty e^{st}d\mathbb{P}(T_{ij}\leq t), \quad i\neq j,~s<0,
\]
we have the well-known result
\begin{equation}\label{Fij}
\hat{F}_{ij}(s) = \frac{\hat{P}_{ij}(s)}{\hat{P}_{jj}(s)}, \quad i\neq j
\end{equation}
(see, for example, \cite[Equation (1.3)]{K57b}).
Karlin and McGregor give in \cite[Equation (3.21)]{K57a} a representation for
$\hat{P}_{ij}(s)$, which upon substitution in \eqref{Fij} yields
\begin{equation} \label{KM}
\hat{F}_{ij}(s) = \frac{Q_i(s)}{Q_j(s)}, \quad 0\leq i<j,
\end{equation}
where $Q_n,~n=0,1,\dots,$ are the {\em birth-death polynomials} associated
with the process $\mathcal{X}$, that is, the $Q_n$ satisfy the recurrence
relation
\begin{equation}
\label{recQ}
\begin{array}{@{}l}
\la_n Q_{n+1}(x) = (\la_n+\mu_n-x)Q_n(x) - \mu_n Q_{n-1}(x), \quad n > 0,\\
\la_0 Q_1(x) = \la_0 + \mu_0 - x, \quad Q_0(x) = 1.
\end{array}
\end{equation}
The representation \eqref{KM} was observed explicitly for the first time by
Karlin and McGregor themselves in \cite[Page 378]{K59}. Since then several authors
have rediscovered the result or provided alternative proofs (see Diaconis and Miclo
\cite{D09} for some references).

In a recent paper in this journal
Gong, Mao and Zhang \cite{G12}, using the theory of Dirichlet forms, extended
Karlin and McGregor's result by establishing a representation for the Laplace
transform of the first-hitting time $T_{ij}$ for {\em any} pair of states $i\neq j$,
as well as asymptotics when either $i$ or $j$ tends to infinity.
It will be shown here that these results may also be obtained by exploiting 
Karlin and McGregor's toolbox, which is the theory of orthogonal polynomials.

Our findings, which are actually somewhat more general than those of Gong, Mao
and Zhang, are presented in Section 3 and proven in Section 4. In the next
section we introduce some further notation, terminology and preliminary
results. Since a path between two states in a birth-death process has to hit
all intermediate states, we obviously have
\[
\hat{F}_{ij}(s) =
\left\{
\begin{array}{l@{}l}
\hat{F}_{0j}(s)/\hat{F}_{0i}(s) &\qquad \mbox{if}~~ i<j\\
\hat{F}_{i0}(s)/\hat{F}_{j0}(s) &\qquad \mbox{if}~~ i>j.
\end{array}
\right.
\]
So for notational simplicity -- and without loss of generality -- we will restrict
ourselves to an analysis of $T_{0n}$ and $T_{n0}$ for $n>0$.

\section{Preliminaries}

We will use the shorthand notation
\[
K_n := \sum_{i=0}^n \pi_i, \quad L_n := \sum_{i=0}^n (\la_i\pi_i)^{-1},
\quad 0\leq n\leq\infty,
\]
and, following Anderson \cite[Chapter 8]{A91},
\begin{equation}
\label{CD}
C := \sum_{n=0}^\infty (\la_n\pi_n)^{-1} K_n, \quad
D := \sum_{n=0}^\infty (\la_n\pi_n)^{-1} (K_\infty - K_n).
\end{equation}
We have $K_\infty+L_\infty=\infty$ by our assumption \eqref{unique}, while,
obviously,
\begin{equation}
\label{CDKL}
K_\infty = \infty ~\lr~ D = \infty, \quad L_\infty = \infty ~\lr~ C = \infty.
\end{equation}
Also, $C+D = K_\infty L_\infty$, so \eqref{unique} is
actually equivalent to $C+D = \infty$. Whether the quantities $C$ and $D$ are
infinite or not determines the type of the boundary at infinity (see, for example,
Anderson \cite[Section 8.1]{A91}), but also, as we shall see, the asymptotic
behaviour of the polynomials $Q_n$ of \eqref{recQ}.

Since the birth-death polynomials $Q_n$ satisfy the three-terms recurrence
relation \eqref{recQ}, they are orthogonal with respect to a positive Borel
measure on the nonnegative real axis, and have {\em positive} and simple
zeros. The orthogonalizing measure for the polynomials $Q_n$ (normalized
to be a probability measure) is not necessarily uniquely determined by the
birth and death rates, but there exists, in any case, a unique {\em natural\/}
measure $\psi$, characterized by the fact that the minimum of its support is
maximal. We refer to Chihara's book \cite{C78} for properties of orthogonal
polynomials in general, and to Karlin and McGregor's papers \cite{K57a} and
\cite{K57b} for results on birth-death polynomials in particular (see also
\cite[Section 3.1]{D15} for a concise overview). For our purposes the
following properties of birth-death polynomials are furthermore relevant.

With $x_{n1}<x_{n2}<\ldots<x_{nn}$ denoting the $n$ zeros of $Q_n(x),$
there is the classical separation result
\[
0 < x_{n+1,i} < x_{ni} < x_{n+1,i+1}, \quad i = 1,2,\ldots,n,~n \geq 1,
\]
so that the limits
\begin{equation}
\label{xi}
\xi_i := \lim_{n\to\infty} x_{ni}, \quad i = 1,2,\ldots.
\end{equation}
exist. We further let
\begin{equation}
\label{si}
\si := \lim_{i\to\infty} \xi_i
\end{equation}
(possibly infinity). The numbers $\xi_i$ may be defined alternatively as
\[
\xi_1 := \inf \mbox{supp}(\psi)
~~\mbox{and}~~
\xi_{i+1} := \inf \{\mbox{supp}(\psi)\cap (\xi_i,\infty)\}, \quad i \geq 1,
\]
where supp stands for {\em support}. So knowledge of the (natural) orthogonalizing
measure for the polynomials $Q_n$ implies knowledge of the numbers $\xi_i$.
It is clear from the definition of $\xi_i$ that
\[
0 \leq \xi_i \leq \xi_{i+1} \leq \si, \quad i \geq 1.
\]
Moreover we have, for all $i \geq 1$,
\[
\xi_{i+1} = \xi_i ~~\llr~~ \xi_i = \si,
\]
as is evident from the alternative definition of $\xi_i$. By suitably interpreting
\cite[Equations (2.6) and (2.11)]{D86} it follows that
\[
\sum_{i=1}^\infty \xi_i^{-1} =
\lim_{n\to\infty}\frac{1}{1+\mu_0L_n}
\sum_{j=0}^n(\la_j\pi_j)^{-1} \sum_{i=0}^j\pi_i(1+\mu_0L_{i-1}),
\]
where the left-hand side should be interpreted as infinity if $\xi_1 = 0$. In particular,
\[
\mu_0 = 0 ~~\lr~~  \sum_{i=1}^\infty\xi_i^{-1} = C.
\]
Also, by \cite[Theorem 2]{D86},
\begin{equation}
\label{sumxi1}
\mu_0 = 0: \quad C < \infty ~\mbox{or}~ 
D < \infty ~\llr~ \sum_{i=2}^\infty\xi_i^{-1} < \infty
\end{equation}
and
\begin{equation}
\label{sumxi}
\mu_0 > 0: \quad C < \infty ~\mbox{or}~ 
D < \infty ~\llr~ \sum_{i=1}^\infty\xi_i^{-1} < \infty.
\end{equation}

Given a sequence of birth-death polynomials $\{Q_n\}$ we obtain the sequence
$\{Q_n^{(l)}\}$ of {\em associated polynomials} of order $l\geq 0$ by
replacing $Q_n$ by $Q_n^{(l)}$, $\la_n$ by $\la_{n+l}$ and $\mu_n$ by
$\mu_{n+l}$ in the recurrence relation \eqref{recQ}.
Evidently, the polynomials $Q_n^{(l)}$ are birth-death polynomials again, so
$Q_n^{(l)}(x)$ has simple, positive zeros
$x_{n1}^{(l)}<x_{n2}^{(l)}<\dots<x_{nn}^{(l)}$ and we can write
\[
Q_n^{(l)}(x) =  
Q_n^{(l)}(0) \prod_{i=1}^n\left(1-\frac{x}{x_{ni}^{(l)}}\right), \quad n,l \geq 0,
\]
while it follows by induction that
\begin{equation}
\label{Qnl0}
Q_n^{(l)}(0) = 1+\mu_l\pi_l(L_{n+l-1}-L_{l-1}), \quad n,l \geq 0,
\end{equation}
where $L_{-1}:=0$. Note that $Q_n^{(0)}(0)=Q_n(0)=1$ for all $n$ if $\mu_0=0$.

Defining the quantities $\xi_i^{(l)}$ and $\si^{(l)}$ in analogy to \eqref{xi} and
\eqref{si}, we have, by \cite[Theorem III.4.2]{C78},
\begin{equation}
\label{xiil}
\xi_i^{(l)} \leq \xi^{(l+1)}_i \leq \xi_{i+1}^{(l)},
\quad l \geq 0,~i\geq 1,
\end{equation}
so that
\[
\si^{(l)} = \si, \quad l \geq 0.
\]
Moreover,
\cite[Theorem 1]{D84} tells us that
\begin{equation}
\label{limxil}
\lim_{l\to\infty}\xi_i^{(l)} = \si, \quad i\geq 1.
\end{equation}

Since the polynomials $Q_n^{(l)}$ are birth-death polynomials they are
orthogonal with respect to a unique {\em natural} (probability) measure $\psi^{(l)}$
on the nonnegative real axis. A key ingredient in our analysis is {\em Markov's
Theorem}, which relates the {\em Stieltjes transform} of the measure $\psi^{(l)}$
to the polynomials $Q_n^{(l)}$ and $Q_n^{(l+1)}$, namely,
\begin{equation}
\label{Markov}
\int_0^\infty \frac{\psi^{(l)}(dx)}{x-s} =
\lim_{n\to\infty} \frac{1}{\la_l}\frac{Q^{(l+1)}_{n-1}(s)}{Q_n^{(l)}(s)},
\quad \mbox{Re}(s)<\xi_1^{(l)}.
\end{equation}
We note that $\psi^{(l)}$ is not necessarily the only orthogonalizing measure for
the polynomials $Q_n^{(l)}$, a setting usually not covered in statements of Markov's
Theorem in the literature (see, for example, \cite[Page 89]{C78}). However, an
extension of the original theorem that serves our needs can be found in Berg \cite{B94}
(see in particular \cite[Section 3]{B94}, where the measure $\mu^{(0)}$ corresponds
to our $\psi^{(l)}$).

We will also have use for a classical result in the theory of continued fractions
relating the Stieltjes transforms of the measures $\psi^{(l)}$ and $\psi^{(l+1)}$,
namely,
\begin{equation}
\label{contfrac}
\int_0^\infty\frac{\psi^{(l)}(dx)}{x-s} =
\left\{\la_l+\mu_l-s-\la_l\mu_{l+1}\int_0^\infty\frac{\psi^{(l+1)}(dx)}{x-s}\right\}^{-1},
~~ \mbox{Re}(s)<\xi_1^{(l)}.
\end{equation}
Again we refer to Berg \cite[Section 4]{B94} for statements of this result in the
generality required in our setting.

Our final preliminary results concern asymptotics for the polynomials $Q_n^{(l)}$ as
$n\to\infty$, which may be obtained by suitably interpreting the results of \cite{K95}
(which extend those of \cite{D86}). We state the results in three propositions and
give more details about their derivations in Section 4. Recall that $\xi^{(l)}_0=-\infty$
and $Q_n(0)=1$ if $\mu_0=0$.  

\begin{proposition}\rm~
\label{KL1}
Let $K_\infty = L_\infty = \infty$. Then $C = D = \infty$, $\si = 0$ and, for $l\geq 0$,
\[
\lim_{n\to\infty} Q_n^{(l)}(x) = \infty 
\quad \mbox{if}~~x<0.
\]
\end{proposition}

\begin{proposition}\rm~
\label{KL2}
Let $K_\infty = \infty$ and $L_\infty < \infty$. Then $D = \infty$ and,\\
(i) for $l\geq 0$,
\[
\lim_{n\to\infty} Q_n^{(l)}(0) = 1 + \mu_l\pi_l(L_\infty - L_{l-1}) < \infty;
\]
(ii) if $C=\infty$, for $l\geq 0$,
\[
\lim_{n\to\infty} Q_n^{(l)}(x) = 
\left\{
\begin{array}{l@{}l}
\infty &\quad \mbox{if}~~x < 0\\
0 &\quad \mbox{if}~~0 < x \leq \xi_k^{(l)} \mbox{~for~some~}k\geq 1;
\end{array}
\right.
\]
(iii) if $C<\infty$, for $l \geq 0$,
\[
\lim_{n\to\infty} \frac{Q_n^{(l)}(x)}{Q_n^{(l)}(0)} =  
\prod_{i=1}^\infty\left(1-\frac{x}{\xi_i^{(l)}}\right), \quad x \in \mathbb{R},
\]
an entire function with simple, positive zeros $\xi_i^{(l)},~i\geq 1$.
\end{proposition}

\begin{proposition}\rm~
\label{KL3}
Let $K_\infty < \infty$ and $L_\infty = \infty$. Then $C = \infty$ and,\\
(i) for $l=0$ and $\mu_0>0$, or $l \geq 1$,
\[
\lim_{n\to\infty} Q_n^{(l)}{(0)} = \infty;
\]
(ii) if $D=\infty$, for $l \geq 0$,
\[
\lim_{n\to\infty} Q_n^{(l)}(x) = 
\left\{
\begin{array}{l@{}l}
\infty &\quad \mbox{if}~~\xi^{(l)}_{2k} < x \leq \xi^{(l)}_{2k+1} \mbox{~for~some~}k\geq 0\\
-\infty &\quad \mbox{if}~~\xi^{(l)}_{2k+1} < x \leq \xi^{(l)}_{2k+2} \mbox{~for~some~}k\geq 0;
\end{array}
\right.
\]
(iii) if $D<\infty$ and $\mu_0=0$,
\[
\lim_{n\to\infty} \frac{Q_n(x)}{L_{n-1}} =
-x K_\infty \prod_{i=1}^\infty \left(1-\frac{x}{\xi_{i+1}}\right), \quad x \in \mathbb{R},
\]
an entire function with simple zeros $\xi_1=0$ and $\xi_{i+1}>0,~i\geq 1$;\\
(iv) if $D<\infty$, for $l=0$ and $\mu_0>0$, or  $l \geq 1$,
\[
\lim_{n\to\infty} \frac{Q_n^{(l)}(x)}{Q_n^{(l)}(0)} =
\prod_{i=1}^\infty\left(1-\frac{x}{\xi_i^{(l)}}\right), \quad x \in \mathbb{R},
\]
an entire function with simple, positive zeros $\xi_i^{(l)},~i\geq 1$.
\end{proposition}

\section{Results}

Representations for $\mathbb{E}[e^{sT_{0n}}\mathbb{I}_{\{T_{0n}<\infty\}}]$ and
$\mathbb{E}[e^{sT_{n0}}\mathbb{I}_{\{T_{n0}<\infty\}}]$ in terms of the
polynomials $Q_n^{(l)}$ are collected in the first theorem.

\begin{theorem}\rm~
\label{thm1}
We have, for $\mu_0 \geq 0$ and $n\geq 1$,
\begin{equation}
\label{T0n}
\mathbb{E}[e^{sT_{0n}}\mathbb{I}_{\{T_{0n}<\infty\}}] = \frac{1}{Q_n(s)},
\quad s < x_{n,1},
\end{equation}
and, if $C+D=\infty$,
\begin{equation}
\label{Tn0}
\mathbb{E}[e^{sT_{n0}}\mathbb{I}_{\{T_{n0}<\infty\}}] =
\frac{\la_0}{\la_n\pi_n}\lim_{N\to\infty} \frac{Q^{(n+1)}_{N-n}(s)}{Q^{(1)}_N(s)},
\quad s<\xi_1^{(1)}.
\end{equation}
\end{theorem}

Note that for $s<0$ we have $\mathbb{E}[e^{sT_{0n}}\mathbb{I}_{\{T_{0n}<\infty\}}]
=\mathbb{E}[e^{sT_{0n}}]$, so the representation \eqref{T0n} reduces to Karlin
and McGregor's result \eqref{KM}. The explicit representation \eqref{Tn0} is new,
but may be obtained by a limiting procedure from Gong, Mao and Zhang
\cite[Corollary 3.6]{G12}, where a finite state space is assumed.

By choosing $s=0$ in \eqref{T0n} and \eqref{Tn0} and using \eqref{Qnl0}, we
obtain expressions for the probabilities $\mathbb{P}(T_{0n}<\infty)$ and
$\mathbb{P}(T_{n0}<\infty)$ that are in accordance with \cite[an unnumbered
formula on page 387 and Theorem 10]{K57b}. For convenience we state the
results as a corollary of Theorem \ref{thm1}, but remark that a proof of \eqref{PTn0}
on the basis of \eqref{Tn0} would require additional motivation in the case
$\xi_1^{(1)} = 0$. 

\begin{corollary}\label{cor1}\rm (\cite{K57b})~
We have, for $\mu_0 \geq 0$ and $n\geq 1$,
\begin{equation}
\label{PT0n}
P(T_{0n} < \infty) = \frac{1}{1 + \mu_0L_{n-1}},
\end{equation}
and, if $C+D=\infty$,
\begin{equation}
\label{PTn0}
\mathbb{P}(T_{n0}<\infty) = 1-\frac{L_{n-1}}{L_\infty}.
\end{equation}
\end{corollary}

\medskip
After a little algebra \eqref{Tn0} and \eqref{Qnl0} lead to
\[
\mathbb{E}[e^{sT_{n0}}\mathbb{I}_{\{T_{n0}<\infty\}}]  
= \ds \left(1 - \frac{L_{n-1}}{L_\infty}\right) \lim_{N\to\infty}
\ds\frac{Q^{(n+1)}_{N-n}(s)/Q^{(n+1)}_{N-n}(0)}{Q^{(1)}_N(s)/Q^{(1)}_N(0)},
\quad s<\xi_1^{(1)}.
\]
Subsequently applying Propositions \ref{KL1}, \ref{KL2} (iii) and \ref{KL3} (iv) we
obtain the second corollary of Theorem \ref{thm1}.

\begin{corollary}\label{cor2}\rm~
If $C+D=\infty$, but $C<\infty$ or $D<\infty$, then, for $n\geq 1$,
\begin{equation}
\label{Tn02}
\mathbb{E}[e^{sT_{n0}}\mathbb{I}_{\{T_{n0}<\infty\}}] =
\left(1 - \frac{L_{n-1}}{L_\infty}\right)
\frac{\ds\prod_{i=1}^\infty\left(1 - \frac{s}{\xi^{(n+1)}_i}\right)}
{\ds\prod_{i=1}^\infty\left(1 - \frac{s}{\xi^{(1)}_i}\right)},
\quad s<\xi_1^{(1)},
\end{equation}
where the infinite products are entire functions with simple, positive zeros
$\xi^{(n+1)}_i$ and $\xi^{(1)}_i,~i\geq 1$.
\end{corollary}

Assuming a denumerable state space, but under the condition $C=\infty$ and $D<\infty$,
Guo, Mao and Zhang give in \cite[Theorem 5.5 (a)]{G12} a representation for
$\mathbb{E}[e^{sT_{n0}}],~s<0$, which is encompassed by Corollary \ref{cor2}.
Indeed, in this case we have $L_\infty = \infty$, and hence, by \eqref{PTn0},
$\mathbb{P}(T_{n0}<\infty) = 1$.

Asymptotic results for $\mathbb{E}[e^{sT_{0n}}\mathbb{I}_{\{T_{0n}<\infty\}}]$
and $\mathbb{E}[e^{sT_{n0}}\mathbb{I}_{\{T_{n0}<\infty\}}]$ as $n\to\infty$ are
summarized in the second theorem.

\begin{theorem}\rm~
\label{thm2}
We have, for $\mu_0 \geq 0$ and $s< 0$, 
\begin{equation}
\label{T0inf}
\lim_{n\to\infty} \mathbb{E}[e^{sT_{0n}}\mathbb{I}_{\{T_{0n}<\infty\}}] =
\left\{
\begin{array}{l@{}l}
\ds\frac{1}{1+\mu_0 L_\infty} \prod_{i=1}^\infty \frac{\xi_i}{\xi_i - s}
&\quad \mbox{if}~~ C < \infty,~ D = \infty\\
0  &\quad \mbox{if}~~ C = \infty,
\end{array}
\right.
\end{equation}
and
\begin{equation}
\label{Tinf0}
\lim_{n\to\infty} \mathbb{E}[e^{sT_{n0}}\mathbb{I}_{\{T_{n0}<\infty\}}] =
\left\{
\begin{array}{l@{}l}
0  &\qquad \mbox{if}~~ C < \infty,~ D = \infty\\
\ds\prod_{i=1}^\infty\frac{\xi_i^{(1)}}{\xi_i^{(1)}-s} &
\qquad \mbox{if}~~ C = \infty,~ D < \infty.
\end{array}
\right.
\end{equation}
The infinite products in \eqref{T0inf} and \eqref{Tinf0} are reciprocals of entire
functions with simple, positive zeros $\xi_i$ and $\xi^{(1)}_i,~i\geq 1$,
respectively.
\end{theorem}

By \eqref{PT0n} we have
\[
\lim_{n\to\infty} P(T_{0n} < \infty) = \frac{1}{1 + \mu_0L_\infty},
\]
so \eqref{T0inf} implies
\[
\lim_{n\to\infty} \mathbb{E}[e^{sT_{0n}}\,|\,T_{0n}<\infty] =
\prod_{i=1}^\infty \frac{\xi_i}{\xi_i - s} \qquad \mbox{if}~~ C < \infty,~ D = \infty,
\]
which generalizes \cite[Theorem 4.6]{G12} where $\mu_0=0$ is assumed. (At the
end of \cite[Section 4]{G12} the authors remark that the case $\mu_0>0$ may
be treated in a way analogous to the case $\mu_0=0$, but no explicit result is
given.) If $C=\infty$ and $D<\infty$ we must have $L_\infty=\infty$ and hence, by
\eqref{PTn0}, $\mathbb{P}(T_{n0}<\infty) = 1$. So \eqref{Tinf0} implies
\[
\lim_{n\to\infty} \mathbb{E}[e^{sT_{n0}}] =
\prod_{i=1}^\infty \frac{\xi_i^{(1)}}{\xi_i^{(1)} - s}
\qquad \mbox{if}~~ C = \infty,~ D < \infty,
\]
which is \cite[Theorem 5.5 (b)]{G12}.

\section{Proofs}

\subsection{Proofs of Propositions 1--3}

The conclusions regarding $C$ and $D$ in the Propositions \ref{KL1}, \ref{KL2} and
\ref{KL3} are given already in \eqref{CDKL}, while the statements (i) in
Propositions \ref{KL2} and \ref{KL3} are implied by \eqref{Qnl0}. The other
statements follow from results in \cite{K95}, where two cases -- corresponding in
the setting at hand to $\mu_0=0$ and $\mu_0>0$ -- are considered simultaneously
by means of a duality relation involving polynomials $R_n$ and $R_n^*$. The
asymptotic results for $R_n$ may be translated into asymptotics for $Q_n$ if
$\mu_0=0$, while the results for $R_n^*$, suitably interpreted, give asymptotics for
$Q_n$ if $\mu_0>0$, and for $Q_n^{(l)}$ with $l\geq 1$.
Concretely, the statements in Proposition \ref{KL1}, Proposition \ref{KL2} (ii) and
Proposition \ref{KL3} (ii) regarding the case $x<0$ follow from \cite[Lemma 2.4
and Theorems 3.1 and 3.3]{K95}, while the results for $x>0$ are implied by
\cite[Theorems 2.2, 3.6 and 3.8]{K95}.
Proposition \ref{KL2} (iii) follows from \cite[Theorem 3.1]{K95} for $l=0$ and
$\mu_0=0$, and from \cite[Corollary 3.2]{K95} for $l=0$ and $\mu_0>0$, and for
$l\geq 1$. Proposition \ref{KL3} (iii) is implied by \cite[Theorems 2.2, 3.3 and
3.4 (ii)]{K95}, while Proposition \ref{KL3} (iv) is a consequence of
\cite[Corollary 3.2]{K95}. 
Finally, the fact that $\si=0$ in the setting of Proposition
\ref{KL1} is stated, for example, in \cite[Theorem 2.2 (iv)]{K95}.

\subsection{Proof of Theorem \ref{thm1}}

As observed already, substitution in \eqref{Fij} of Karlin and McGregor's formula for
$\hat{P}_{ij}(s)$ given on \cite[Equation (3.21)]{K57a} leads to \eqref{KM} and hence,
by analytic continuation, to \eqref{T0n}.

To obtain \eqref{Tn0} we note that \cite[Equation (3.21)]{K57b}) also yields
\[
\hat{P}_{10}(s) = -\frac{1}{\la_0} + Q_1(s)\hat{P}_{00}(s) =
\frac{1}{\la_0}\left[(\la_0+\mu_0-s)\hat{P}_{00}(s) - 1\right],
\]
which upon substitution in \eqref{Fij} leads to
\[
\hat{F}_{10}(s) =
\frac{1}{\la_0}\left[\la_0+\mu_0-s-\frac{1}{\hat{P}_{00}(s)}\right], \quad s<0.
\]
Moreover, by Karlin and McGregor's representation formula for the transition
probabilities $P_{ij}(t)$ (see \cite[Section III.6]{K57a}) we have $P_{00}(t) =
\int_0^\infty e^{-xt}\psi(dx)$, where $\psi$ is a (probability) measure with
respect to which the polynomials $Q_n$ are orthogonal. Since the condition
$C+D=\infty$ is equivalent to \eqref{unique}, it ensures that the transition
probabilities are uniquely determined by the birth and death rates, whence
$\psi$ must be the natural measure (see \cite{K57a}). So we have
\begin{equation}
\label{P00}
\hat{P}_{00}(s) = \int_0^\infty\frac{\psi(dx)}{x-s}, \quad s<\xi_1.
\end{equation}
Subsequently applying \eqref{contfrac} with $l=0$, it follows that
\begin{equation}
\label{F10}
\hat{F}_{10}(s) = \mu_1\int_0^\infty\frac{\psi^{(1)}(dx)}{x-s}, \quad s<0,
\end{equation}
whence, more generally,
\[
\hat{F}_{l,l-1}(s) = \mu_l\int_0^\infty\frac{\psi^{(l)}(dx)}{x-s},
\quad s<0,~l\geq 1,
\]
and, by analytic continuation,
\begin{equation}
\label{Ell-1}
\mathbb{E}[e^{sT_{l,l-1}}\mathbb{I}_{\{T_{l,l-1}<\infty\}}] =
\mu_l\int_0^\infty\frac{\psi^{(l)}(dx)}{x-s}, \quad s<\xi_1^{(l)},~l\geq 1.
\end{equation}
Since $T_{n0} = T_{n,n-1}+\dots +T_{10}$, while $T_{n,n-1},\dots,T_{10}$ are
independent random variables, Markov's Theorem \eqref{Markov} implies that we
can write
\[
\begin{array}{@{}l@{}l}
\mathbb{E}[e^{sT_{n0}}\mathbb{I}_{\{T_{n0}<\infty\}}] &~
= \mathbb{E}[e^{sT_{n,n-1}}\mathbb{I}_{\{T_{n,n-1}<\infty\}}] \dots
\mathbb{E}[e^{sT_{10}}\mathbb{I}_{\{T_{10}<\infty\}}]\\
&~= \ds\frac{\mu_1\dots\mu_n}{\la_1\dots\la_n}
\ds\left(\lim_{N\to\infty}\frac{Q^{(n+1)}_{N-n}(s)}{Q^{(n)}_{N-n+1}(s)}\right)\dots
\ds\left(\lim_{N\to\infty}\frac{Q^{(2)}_{N-1}(s)}{Q^{(1)}_N(s)}\right)\\
&~= \ds\frac{\la_0}{\la_n\pi_n}\lim_{N\to\infty} \frac{Q^{(n+1)}_{N-n}(s)}{Q^{(1)}_N(s)}.
\end{array}
\]
Recalling \eqref{xiil} we conclude that this expression holds for $s<\xi_1^{(1)}$.

\subsection{Proof of Theorem \ref{thm2}}

Letting $n\to\infty$ in \eqref{T0n} and applying the results of Propositions
\ref{KL1}, \ref{KL2} and \ref{KL3} readily yields the first statement of
Theorem \ref{thm2}.

To prove the second statement we employ Corollary \ref{cor2}. First note that, for
$a>0$ and $s\leq 0$, we have $1 \leq 1-\frac{s}{a}\leq e^{-s/a}$, so that, for $l \geq 0$,
\[
1 \leq \prod_{i=1}^\infty \left(1 - \frac{s}{\xi^{(l)}_i}\right)
\leq \exp\left\{-s\sum_{i=1}^\infty \frac{1}{\xi^{(l)}_i}\right\}, \quad s \leq 0,
\]
provided $\xi_1^{(l)}>0$.
Defining $C^{(l)}$ and $D^{(l)}$ in analogy to \eqref{CD} it is easily seen that
\[
C < \infty ~\llr~ C^{(l)} < \infty, \quad D < \infty ~\llr~ D^{(l)} < \infty.
\]
So, assuming $C<\infty$ or $D<\infty$, we have, by \eqref{sumxi},
\[
\sum_{i=1}^\infty \frac{1}{\xi^{(l)}_i} < \infty, \quad l\geq 1.
\]
Hence $\si^{(l)} = \si =\infty$, so that, by \eqref{limxil}, $\xi_i^{(l)} \to \infty$
as $l \to \infty$. As a consequence
\[
\lim_{l\to\infty} \ds\prod_{i=1}^\infty\left(1 - \frac{s}{\xi^{(l)}_i}\right) = 1,
\quad s \leq 0,
\]
and the result follows since $L_\infty<\infty$ if $C<\infty$, whereas
$L_\infty = \infty$ if $C=\infty$ and $D<\infty$.

\section{Concluding remarks}

First we note that the result \eqref{F10} -- or rather a generalization of
\eqref{F10} -- may be derived directly from the Kolmogorov differential equations and
Karlin and McGregor's representation formula for the transition probabilities $P_{ij}(t)$.
The argument is given on \cite[Page 508]{D03b} (and essentially already on
\cite[Page 385]{K57b}) and yields
\[
\mathbb{P}(t < T_{n0} < \infty) =
\mu_1 \int_0^\infty \frac{e^{-xt}}{x}Q^{(1)}_{n-1}(x)\psi^{(1)}(dx), \quad n\geq 1,
\]
so that
\begin{equation}\label{Fn0}
\hat{F}_{n0}(s) = \mu_1 \int_0^\infty \frac{Q^{(1)}_{n-1}(x)}{x-s}\psi^{(1)}(dx),
\quad s < 0,~n\geq 1.
\end{equation}
Note that as a consequence of \eqref{Tn0} and \eqref{Fn0} we have, for all $m \geq 0$,
\[
\int_0^\infty \frac{Q^{(1)}_m(x)}{x-s}\psi^{(1)}(dx) = 
\frac{\pi_1}{\la_{m+1}\pi_{m+1}}\lim_{n\to\infty} \frac{Q^{(m+2)}_{n-m-1}(s)}{Q^{(1)}_n(s)},
\quad s < 0,
\]
which implies a partial extension of Markov's Theorem \eqref{Markov} to the effect that,
for $m \geq 0$ and $l \geq 1$ (and $l=0$ if $\mu_0>0$),
\begin{equation}
\label{Markovgen}
\int_0^\infty \frac{Q^{(l)}_m(x)}{x-s}\psi^{(l)}(dx) = 
\frac{\pi_l}{\la_{m+l}\pi_{m+l}}\lim_{n\to\infty} \frac{Q^{(m+l+1)}_{n-m-1}(s)}{Q^{(l)}_n(s)}.
\end{equation}
Using \eqref{Markov}, \eqref{contfrac}, and the recurrence relation for the
polynomials $Q^{(l)}_n$ it may be shown by induction that \eqref{Markovgen} is actually
valid for all $l \geq 0$, $m \geq 0$ and Re$(s) < \xi^{(l)}_1$. Substitution $s=0$ in
\eqref{Markovgen} and \eqref{Qnl0} lead in particular to
\[
\int_0^\infty \frac{Q_m(x)}{x}\psi(dx) = 
\frac{1}{\la_m\pi_m}\lim_{n\to\infty} \frac{Q^{(m+1)}_{n-m-1}(0)}{Q_n(0)} =
\frac{L_\infty-L_{m-1}}{1+\mu_0L_\infty},
\]
which is consistent with \cite[Equations (9.9) and (9.14)]{K57b} (also when $\xi_1=0$).

If we do not impose the condition $C+D=\infty$, the birth and death rates do
not necessarily determine a birth-death process uniquely. However, as observed
in \cite{G12}, several results remain valid if $C+D<\infty$, provided they are
interpreted as properties of the {\em minimal} process, which is the process with
an absorbing boundary at infinity (and which is always associated with the {\em
natural} measure for the polynomials $Q_n$, see \cite{D87}). Concretely, if
$C+D<\infty$ the arguments leading to Theorem \ref{thm1}, and hence Theorem
\ref{thm1} itself and Corollary \ref{cor1}, remain valid. Moreover, the results in
\cite{K95} imply that, for $l\geq 0$,
\[
C+D<\infty ~\lr~ \lim_{n\to\infty} \frac{Q_n^{(l)}(x)}{Q_n^{(l)}(0)} =
\prod_{i=1}^\infty\left(1-\frac{x}{\xi_i^{(l)}}\right),
\]
an entire function with simple, positive zeros
$\xi_i^{(l)},~i\geq 1$. (Note that this complements Propositions \ref{KL1} --
\ref{KL3}.) Hence also \eqref{Tn02} remains valid. Finally, letting $n\to\infty$
in \eqref{T0n} and \eqref{Tn02} we readily conclude that the results in Theorem
\ref{thm2} for $C<\infty,~D=\infty$ are actually valid for $C+D<\infty$
as well.
	
In the setting $C+D<\infty$ Gong, Mao and Zhang \cite{G12} pay attention also
to the {\em maximal\/} process, the process that is characterized by a reflecting
barrier at infinity. In this case the measure featuring in the representation
for $P_{00}(t)$, and hence in \eqref{P00}, is {\em not\/} the natural measure.
Although, applying the results of \cite{D87}, the relevant measure can be identified
and expressed in terms of a {\em natural} measure corresponding to a {\em dual\/}
birth-death process, application of Markov's Theorem does not seem feasible in this
case.

Our final remark is the following. Choosing $l=0$, letting $s \uparrow \xi_1$
in \eqref{contfrac}, and using the recurrence relation \eqref{recQ}, we readily get
\[
\int_0^\infty \frac{\psi(dx)}{x-\xi_1} = \infty ~\llr~
\mu_1 \int_0^\infty \frac{\psi^{(1)}(dx)}{x-\xi_1} = Q_1(\xi_1).
\]
In fact, using \eqref{contfrac} again, it is not difficult to generalize this result to
\[
\int_0^\infty \frac{\psi(dx)}{x-\xi_1} = \infty ~\llr~
\mu_l \int_0^\infty \frac{\psi^{(l)}(dx)}{x-\xi_1} =
\frac{Q_l(\xi_1)}{Q_{l-1}(\xi_1)}, \quad l\geq 1,
\]
which upon substitution in \eqref{Ell-1} leads to
\[
\int_0^\infty \frac{\psi(dx)}{x-\xi_1} = \infty ~\lr~
\mathbb{E}[e^{\xi_1 T_{n0}}\mathbb{I}_{\{T_{n0}<\infty\}}] = Q_n(\xi_1), \quad n\geq 1.
\]
(Since, by \cite[Theorem 3.1]{D06}, the condition above is equivalent to
$\xi_1$-recurrence of the process, this result may also be obtained by applying
\cite[Lemma 3.3.3 (iii)]{J95} to the setting at hand, see \cite[Lemma 3.2]{G15}.)
It now follows from \eqref{Tinf0} that
\[
C = \infty,~ D < \infty ~\lr~
\lim_{n\to\infty} Q_n(\xi_1) =
\ds\prod_{i=1}^\infty\frac{\xi_i^{(1)}}{\xi_i^{(1)} - \xi_1} < \infty.
\]
If $\mu_0=0$ then $\xi_1=0$, so the result does not take us by surprise, but for
$\mu_0>0$ we regain an interesting extension of Proposition \ref{KL3} (iv) -- recently
obtained by Gao and Mao \cite[Lemma 3.4]{G15} -- since it has consequences for
the existence of quasi-stationary distributions.


\end{document}